\documentclass[12pt,reqno]{amsart}

\usepackage[margin=1in]{geometry}  
\usepackage{amsmath}
\usepackage{amssymb}
\usepackage{verbatim}   
\usepackage{color}    
\usepackage{mathrsfs}  
\usepackage[colorlinks=true]{hyperref}
\hypersetup{urlcolor=blue, citecolor=blue, pdfstartview=FitH}
\hypersetup{pdfstartview={XYZ null null 1.25}}
\usepackage{fontenc}
\usepackage{textcomp}
\usepackage{fix-cm}
\usepackage{array}
\usepackage{enumerate}
\usepackage{latexsym, amsfonts, amssymb}
\usepackage{color}

\numberwithin{equation}{section}

\newtheorem{theorem}[equation]{Theorem}

\newtheorem{remark}[equation]{Remark}

\begin{document}

\title[Damped Euler Equation with Three-Point Boundary Condition]{Blowup of Solutions to a Damped Euler Equation with Homogeneous Three-Point Boundary Condition}

\author{Ikechukwu Obi-Okoye}
\address{Ikechukwu Obi-Okoye\newline Department of Physics\newline University of North Georgia\newline 82 College Circle\newline Dahlonega, GA  30597}
\email{ivobio3570@ung.edu}

\author{Alejandro Sarria}
\address{Alejandro Sarria\newline Department of Mathematics and Statistics\newline University of North Georgia\newline 82 College Circle\newline 
	Dahlonega, GA  30597}
\email{alejandro.sarria@ung.edu}

\begin{abstract}
In \cite{Yuen}, it was established that solutions to the inviscid Proudman-Johnson equation subject to a homogeneous three-point boundary condition can develop singularities in finite time. In this paper, we consider the possibility of singularity formation in solutions of the generalized, inviscid Proudman-Johnson equation with damping subject to the same homogeneous three-point boundary condition. In particular, we derive conditions the initial data must satisfy in order for solutions to blowup in finite time with either bounded or unbounded smooth damping term.   
\end{abstract}

%
%
%
%


\maketitle

\section{Introduction}
\label{sec:intro}
We study the initial value problem (IVP)
\begin{equation}
\label{3p}
\begin{cases}
u_{xxt}+uu_{xxx}+\beta u_xu_{xx}+\alpha(t)u_{xx}=0,\,\,\,\,\,\,\,\,\,&x\in[0,1],\,\,t>0,
\\
u(x,0)=u_0(x),\,\,\,\,\,\,&x\in[0,1],
\end{cases}
\end{equation}
where $\beta=\frac{n-3}{n-1}$, $n\in\mathbb{Z}^+$, $n\geq2$, and $u(x,t)$ satisfies the homogeneous three-point boundary condition
\begin{equation}
\label{bc}
u(1,t)=u_x(0,t)=u_x(1,t)=0
\end{equation}

Equation (\ref{3p})i can be obtained by imposing on the incompressible $n-$dimensional Euler equations with damping
$$\boldsymbol{u}_t+(\boldsymbol{u}\cdot\nabla)\boldsymbol{u}+\alpha(t)\boldsymbol{u}=-\nabla p,\qquad\quad\nabla\cdot \boldsymbol{u}=0$$ 
velocities of the form $$\boldsymbol{u}(x,\boldsymbol{x}^\prime,t)=(u(x,t),-\frac{\boldsymbol{x}^\prime}{n-1}u_x(x,t))$$

for $\boldsymbol{x}^\prime=\{x_2,...,x_n\},$ or by using the cylindrical coordinate representation 

$$u^r=-\frac{r}{n-1}\,u(x,t),\quad u^{\theta}=0,\qquad u^x=u(x,t)$$

where $r=\left|\boldsymbol{x}^\prime\right|$ (\cite{Childress}, \cite{Weyl1}, \cite{Saxton1}, \cite{Okamoto1}, \cite{Escher1}, \cite{Proudman1}).

We will refer to (\ref{3p})-(\ref{bc}) as the initial boundary value problem for the generalized, inviscid Proudman-Johnson equation with damping. 

We remark that the undamped version of  (\ref{3p})i), i.e.
\begin{equation}
\label{undamped}
u_{xxt}+uu_{xxx}+\beta u_xu_{xx}=0
\end{equation}
arises in several physical and geometrical contexts. For instance, when $\beta=3,$ (\ref{undamped}) reduces to the Burgers' equation of gas dynamics. For $\beta=2,$ it becomes the Hunter-Saxton equation (HS) describing the orientation of waves in massive nematic liquid crystals (\cite{Hunter1}, \cite{Bressan1}, \cite{Dafermos1}, \cite{Yin1}). From a more geometric point of view, periodic solutions to the HS equation also describe geodesics on the group $\mathcal{D}(\mathbb{S})\backslash Rot(\mathbb{S})$ of orientation preserving diffeomorphisms on the unit circle modulo rigid rotations with respect to the right-invariant metric (\cite{Khesin1}, \cite{Bressan1}, \cite{Tiglay1}, \cite{Lenells1})
$$\left\langle f,g\right\rangle=\int_{\mathbb{S}}{f_xg_xdx}$$

Furthermore, Lenells and Misiolek \cite{Lenells2} showed that, for any $\beta$, (\ref{undamped}) arises as the geodesic equation of the affine connection $\nabla^{(\alpha)}$ on the group $\mathcal{D}(\mathbb{S})\backslash Rot(\mathbb{S})$. See also \cite{Bauer1} for yet another derivation of (\ref{undamped}) as a geodesic equation.   

From a more heuristic point of view, (\ref{3p})i) may serve as a tool to better understand the role that convection and stretching play in the regularity of solutions to one-dimensional fluid evolution equations; it has been argued that the convection term can sometimes cancel some of the nonlinear effects and contribute positively to regularity of solutions (\cite{Ohkitani1}, \cite{Hou1}, \cite{Okamoto3}). More particularly, setting $\omega=u_{xx}$ in (\ref{3p})i) yields
\begin{equation}
\label{eq:diff}
\omega_{t}+\underbrace{u\omega_{x}}_{\text{convection}}+\beta\underbrace{\omega u_{x}}_{\text{stretching}}+\alpha(t)\omega=0
\end{equation}

In the undamped case $(\alpha(t)\equiv 0)$ and $\beta=-1$, (\ref{eq:diff})  becomes a one-dimensional analogue of the three-dimensional vorticity equation of incompressible inviscid fluids
$$\boldsymbol{\omega}_t+\boldsymbol{u}\cdot\nabla\boldsymbol{\omega}=\boldsymbol{\omega}\cdot\nabla\boldsymbol{u},\qquad \boldsymbol{\omega}=\nabla\times\boldsymbol{u}$$

The nonlinear terms in equation (\ref{eq:diff}) represent the competition between nonlinear convection and stretching (\cite{Holm1}, \cite{Wunsch5}). More particularly, the parameter $\beta\in\mathbb{R}$ is related to the ratio of stretching to convection.

\section{Previous Regularity Results}

The global regularity of solutions of the damped equation (\ref{3p})i) for a particular value of $\beta$ is discussed in \cite{Sarria3} for solutions satisfying periodic boundary conditions. In the undamped case, finite-time singularity formation has been extensively studied for periodic solutions as well as Dirichlet boundary conditions (see e.g. \cite{Sarria1}, \cite{Sarria2} and references therein). However, we are unaware of any results concerning singularity formation in solutions of the damped equation (\ref{3p})i) with the homogeneous three-point boundary condition (\ref{bc}). In the undamped case with homogeneous three-point boundary condition (\ref{bc}), the only result we are aware of is Theorem (\ref{thm:Yuen}) below established by Yuen (\cite{Yuen}) where conditions on the initial data leading to finite-time blowup are derived for $\beta=-1$. Therefore, it is of interest to study whether finite-time blowup of solutions is still a possibility after incorporating damping into the system as well as determining if varying the value of the parameter $\beta$ has any effect on the regularity of solutions which, in turn, could lead to a better understanding of the competing effects between nonlinear convection and stretching as discussed at the end of the previous section.

\begin{theorem}
\label{thm:Yuen} 
(Yuen \cite{Yuen}) Consider the $C^3$ solutions to the IBVP (\ref{3p})-(\ref{bc}) for the inviscid undamped Proudman Johnson equation ($n=2,\, \alpha(t)\equiv0$). If the initial velocity satisfies
$$U_0=\int_0^1u_0'(x)\,dx=-u_0(0)>0$$

then the solutions blowup on or before the finite time $1/(2U_0)$.
\end{theorem}

The outline of the paper is a follows. In Section \ref{bounded}, we prove finite-time blowup of solutions to (\ref{3p})-(\ref{bc}) with smooth and bounded time-dependent damping term, while the formation of singularities in finite time with smooth and unbounded damping term is discussed in Section \ref{unbounded}. 

\section{Finite-time Blowup with Smooth, bounded Damping}
\label{bounded}

In this section, we will derive conditions on the initial data $u_0(x)$ which show that the presence of a smooth (bounded or unbounded) time-dependent damping term $\alpha(t)$ is insufficient to arrest finite-time blowup for any value of the parameter $\beta\in\mathbb{R}$.

\begin{theorem}
\label{thm:ike} 
Consider the IBVP (\ref{3p})-(\ref{bc}) for the generalized inviscid damped Proudman Johnson equation with initial data $u_0(x)\in C^{\infty}([0,1])$ and smooth bounded damping term $\alpha(t):[0,\infty)\to\mathbb{R}^+$.

Set  
$$H_0=-\int_0^1u_0'(x)\,dx=u_0(0)$$

and
\begin{equation}
\label{t2}
M:=\sup_{t\in[0,\infty)}\alpha(t)
\end{equation}

for some $M\in\mathbb{R}^+$. If 
\begin{equation}
\label{t1}
H_0<0\quad\text{and, more particularly,}\quad H_0<M\left(\frac{1-n}{n}\right)
\end{equation}

for fixed $n\in\mathbb{Z}^+$, $n\geq2$, then 

\begin{equation}
\label{t3}
\lim_{t\nearrow t^*}u(0,t)=-\infty
\end{equation}

for positive $t^*$ given by 
$$t^*=-\frac{1}{M}\ln\left(1-\frac{M(1-n)}{nH_0}\right)$$

%
%
%
%

\begin{proof}

Multiplying (\ref{3p})i) by $x$, integrating over $x\in[0,1]$ and using (\ref{bc}) gives
\begin{equation}
\label{eq1}
H'(t)+\alpha(t)H(t)+\frac{n}{n-1}\|u_x(\cdot,t)\|^2_{2}=0
\end{equation}
for
\begin{equation}
\label{H}
H(t)=-\int_0^1u_x(x,t)\,dx=u(0,t)
\end{equation}

Now, the Cauchy-Schwarz inequality implies that 
$$H^2\leq\|u_x(\cdot,t)\|^2_{2}$$

and so equation (\ref{eq1}) yields the inequality
\begin{equation}
\label{eq2}
H'(t)+\alpha(t)H(t)+\frac{n}{n-1}H^2\leq 0
\end{equation}

Setting $f=H^{-1}$ in (\ref{eq2}) and using an integrating factor argument we obtain
$$\frac{d}{dt}\left(f(t)e^{-\int_0^t\alpha(s)\,ds}\right)\geq\frac{n}{n-1}e^{-\int_0^t\alpha(s)\,ds}$$ 

which gives, after integrating,
\begin{equation}
\label{eq3}
\frac{1}{H(t)}\geq\frac{n(1+H_0g(t))}{(n-1)H_0g'(t)}
\end{equation}

for $H_0=H(0)$ and
\begin{equation}
\label{eq10}
g(t)=\frac{n}{n-1}\int_0^te^{-\int_0^s\alpha(z)\,dz}\,ds
\end{equation}
 
Next, from (\ref{t2}), 
$$\alpha(t)\leq M$$

for all $t\in[0,\infty)$ for some $M\in\mathbb{R}^+$. Consequently,
$$g'(t)\geq \frac{n}{n-1}e^{-Mt}$$
and, after integrating,
$$g(t)\geq\frac{n}{n-1}\left(\frac{1-e^{-Mt}}{M}\right)$$

Using these inequalities on (\ref{eq3}) we obtain
\begin{equation}
\label{eq4}
\frac{1}{H(t)}\geq\Lambda(t)
\end{equation}

for
\begin{equation}
\label{Lambda}
\Lambda(t)=\frac{nN(t)}{(n-1)^2MH_0g'(t)}
\end{equation}
and
$$N(t)=M(n-1)+nH_0(1-e^{-Mt})$$

Now, (\ref{t1})ii) implies that
$$\lim_{t\to\infty}N(t)=M(n-1)+nH_0<0$$

This, along with 
$$N(0)=M(n-1)>0,\quad N'(t)=nMH_0e^{-Mt}<0\quad \text{and}\quad N(t)\in C([0,\infty))$$

implies the existence of a finite time $t^*>0$ such that
\begin{equation}
\label{eq5}
\lim_{t\nearrow t^*}N(t)=0
\end{equation}

In turn, the above argument implies that $\Lambda(t)<0$ when $t\in[0,t^*)$ and
$$\lim_{t\nearrow t^*}\Lambda(t)=0$$

Finally, (\ref{eq2}) gives
$$H'(t)+\alpha(t)H(t)\leq 0$$

which yields 
\begin{equation}
\label{eq6}
H(t)\leq H_0\,e^{-\int_0^t\alpha(s)ds}
\end{equation}

Since $H_0<0$, (\ref{eq6}) implies that $H(t)<0$ for as long as it exists. Thus
\begin{equation}
\label{eq7}
0\geq\frac{1}{H(t)}\geq \Lambda(t)
\end{equation}

and the Theorem follows by letting $t$ approach $t^*$ in (\ref{eq7}).

%
%

\end{proof}

\end{theorem}

\section{Finite-time Blowup with Smooth, unbounded Damping}
\label{unbounded}

In Theorem \ref{thm:ike} of the previous section, finite-time blowup was established for smooth and bounded time-dependent damping. In Theorem \ref{thm:ike2} below we show that finite-time blowup is still possible with smooth but unbounded time-dependent damping. 

\begin{theorem}
\label{thm:ike2} 
Consider the IBVP (\ref{3p})-(\ref{bc}) for the generalized inviscid damped Proudman Johnson equation with initial data $u_0(x)\in C^{\infty}([0,1])$. Let $E_1(x)$ be the exponential integral

\begin{equation}
\label{eq12}
E_1(x)=\int_x^{\infty}\frac{e^{-t}}{t}dt,\qquad x>0
\end{equation}

If

\begin{equation}
\label{eq8}
u_0(0)<-\frac{c(n-1)}{ne^{1/c}E_1(1/c)}
\end{equation}

for $c=\alpha'(0)\in\mathbb{R}^+$ and $n\in\mathbb{Z}^+$, $n\geq 2$, then there exists smooth, unbounded damping $\alpha(t):[0,\infty)\to\mathbb{R}^+$ and a finite time $t^*>0$ such that

\begin{equation}
\label{eq13}
\lim_{t\nearrow t^*}u(0,t)=-\infty
\end{equation}

\begin{proof}
	
Suppose   
$$H_0=-\int_0^1u_0'(x)\,dx=u_0(0)<0$$

From (\ref{eq3}) and (\ref{eq6}) we have that

\begin{equation}
\label{eq9}
0>\frac{1}{H(t)}\geq \frac{n(1+H_0g(t))}{(n-1)H_0g'(t)}
\end{equation}
	
for $g(t)$ as in (\ref{eq10}) and $n\in\mathbb{Z}^+$, $n\geq 2$. Set
$$\alpha(t)=e^{ct}$$

for some $c\in\mathbb{R}^+$. Then, after some straightforward computations, we can write (\ref{eq9}) as
 	
\begin{equation}
\label{eq11}
0>\frac{1}{H(t)}\geq \frac{1+H_0\phi(n)\eta(t)}{e^{\frac{1}{c}}H_0e^{-\frac{e^{ct}}{c}}}
\end{equation}

for
$$\phi(n)=\frac{ne^{\frac{1}{c}}}{c(n-1)}$$

and

$$\eta(t)=\int_{\frac{1}{c}}^{\frac{e^{ct}}{c}}\frac{e^{-u}}{u}\,du$$

Note that $\eta(0)=0$, while $\eta'(t)>0$ and $\eta''(t)<0$ for all $t\in[0,\infty)$. Moreover, for the exponential integral $E_1(x)$ as defined in (\ref{eq12}), we have that
$$\lim_{t\to\infty}\eta(t)=E_1(1/c)$$

The above remarks, along with the fact that the denominator of the right hand-side of (\ref{eq11}) is negative for all $t\in[0,\infty)$, imply that if $H_0=u_0(0)$ satisfies (\ref{eq8}), then there exists a finite time $t^*>0$ such that (\ref{eq13}) follows. 

\end{proof}	
\end{theorem}

\begin{remark}
Below are some sample upper bounds (\ref{eq8}) (computed using Mathematica) for the initial data $H_0$ in Theorem \ref{thm:ike2} for certain values of $n$ and $c$.

\begin{itemize}
\item If $n=2$ and $c=1$, then 

$$H_0<-\frac{1}{2eE_1(1)}\approx-0.84$$

\item If $n=2$ and $c=0.01$, then 

$$H_0<-\frac{0.01}{2e^{100}E_1(100)}\approx-0.51$$

\end{itemize}

\end{remark}

\section{Conclusions and Open Questions}

In this paper, we derived conditions on the initial data which guarantee the formation of singularities in finite time in solutions of the generalized inviscid Proudman-Johnson when smooth and bounded damping is incorporated into the system. Similarly, we have shown that finite-time blowup is also possible if the damping term is unbounded. It is unknown if solutions will persist for all time (instead of blowing up in finite time) if the initial data satisfies the opposite of inequality (\ref{t1})ii) in Theorem \ref{thm:ike}, namely, if
$$H_0\geq M\left(\frac{1-n}{n}\right)$$

It is also of interest to investigate how solutions to other fluid-related models such as the Boussinesq equations or the MHD equations behave under a homogeneous three-point boundary condition.

\end{document}